\documentclass[reqno,a4paper,11pt]{amsart}
\usepackage{amssymb}
\usepackage[dvips, lmargin=4cm, rmargin=4cm, tmargin=4cm, bmargin=4cm]{geometry}
\usepackage[colorlinks=true, citecolor=red, linkcolor=blue, urlcolor=blue]{hyperref}
\usepackage{mathrsfs}
\usepackage{cite}
\usepackage{xfrac}
\usepackage{graphicx}
\usepackage{hyperref}
\usepackage{amsmath}
\usepackage[T1]{fontenc}
\usepackage{enumitem}
\usepackage{dsfont}
\everymath{\displaystyle}
\newtheorem{theorem}{Theorem}[section]
\newtheorem{definition}{Definition}[section]
\newtheorem{remark}{Remark}[section]
\newtheorem{lemma}{Lemma}[section]
\newtheorem{proposition}{Proposition}[section]

\numberwithin{equation}{section}

\begin{document}
\title[Fractional Caputo-derivative Orlicz  space]{On a novel  fractional Caputo-derivative Orlicz space}
\author[Kasmi Ayoub, Azroul El Houssine, and  Shimi Mohammed]
{Kasmi Ayoub $^1$, Azroul El Houssine  $^2$,   and  Shimi Mohammed $^3$ }
\address{Ayoub Kasmi, Elhoussine Azroul \newline
 Sidi Mohamed Ben Abdellah
 University,
 Faculty of Sciences Dhar el Mahraz, Laboratory of Mathematical Analysis and Applications, Fez, Morocco.}
 \email{$^1$ayoub.kasmi@usmba.ac.ma}
 \email{$^2$elhoussine.azroul@usmba.ac.ma}
\address{Mohammed Shimi\newline
 Sidi Mohamed Ben Abdellah
 University,
ENS of Fez, Laboratory of Mathematical Analysis and Applications, Fez, Morocco.}
\email{$^3$mohammed.shimi2@usmba.ac.ma}
\date{Month, Day, Year}
\begin{abstract}
In this work, we aim to explore whether a novel type of fractional space can be defined using Orlicz spaces and fractional calculus. This inquiry is fruitful, as extending classical results to new contexts can lead to a better and deeper understanding of those classical results. Our main objective is to introduce a new fractional Caputo-derivative Orlicz space, denoted by $\mathcal{O}^{\alpha,G}(\Lambda, \mathbb{R}) $. We are interested in several qualitative properties of this space, such as reflexivity, completeness, and separability. Additionally, we establish a continuous embedding results of this space into a suitable Orlicz space and the space of continuous functions. As an application, we use the mountain-pass theorem (MPT) to ensure the existence of a nontrivial weak solution for a new class of fractional-type problems in Orlicz space.      

\end{abstract}
\maketitle
\tableofcontents
\section{Introduction}\label{sec1}
Fractional calculus constitutes a significant research area because it extends traditional calculus by allowing the use of non-integer order derivatives and integrals, providing powerful tools for modelling and solving complex problems across various fields including applied mathematics, physics, biology, engineering, economics, and others (see for instance \cite{Freed,Lundstrom,Glockle,Hilfer,Mainardi, Kirchner}). \\
On other hand, the mathematical analysis of complex phenomena relies on suitable functional spaces with specific properties. Among these spaces,  Orlicz-Sobolev  spaces \cite{Adams,MusiElaK,Rutickii,Mihailescu,Rao} provide a powerful framework for extending the analysis beyond classical Lebesgue and Sobolev spaces, allowing for a wider range of functions and growth behaviors. Their main advantage lies in their ability to capture the finer properties of functions, where traditional Sobolev spaces sometimes struggle to represent them optimally. In addition, Orlicz-Sobolev spaces are easier to handle analytically because of their convexity and modularity properties, making them easier to handle mathematically. \\

\noindent Alongside, fractional Sobolev spaces,  introduced by  Nikol'skij \cite{Nik1} in 1951,  are an extension of classical Sobolev spaces to accommodate fractional (non-integer) orders of differentiability offering a bridge between classical analysis and more complex, non-local phenomena. For $s\in (0,1)$,  $p\in(1,+\infty),$ these spaces are defined as
$$
W^{s, p}(\Omega)=\left\{u \in L^{p}(\Omega): \tfrac{u(x)-u(y)}{|x-y|^{\tfrac{{\tiny N}}{p}+s}} \in L^{p}(\Omega \times \Omega)\right\},
$$
where $\Omega \subset \mathbb{R}^{N}$ is an open set. They  are useful in various fields such as partial differential equations, image processing, and fractal geometry. \\
In the last few years, based on the same approach given in aforementioned spaces, some scholars have extended the Orlicz-Sobolev spaces to include tha fractional case, see for example \cite{Bonder,Shimi1,Srati} and the references therein.  \\
% The combination of these two fields  Orlicz spaces and fractional calculus opens the door to a new space, which we will henceforth call \textbf{"Orlicz fractional derivative space."}\\

In \cite{Jiao} Jiao and Zhou proposed an other way to define the Sobolev space in the fractional case. They introduced the space $E_0^{\alpha, p}$ as a set of functions $v \in$ $L^p\left(\Lambda, \mathbb{R}^N\right)$ having an
$\alpha$-order Caputo fractional derivative  ${ }_0^C D_t^\alpha v$ in $L^p\left(\Lambda, \mathbb{R}^N\right)$ and $v(0)=v(T)=0$ with $\Lambda = [0,T]$. More recently, José el al. \cite{Costa} introduced a new space using a novel derivative called "$\psi$-Hilfer fractional derivative ${ }^{\mathbf{H}} \mathbf{D}_{0+}^{\alpha, \beta; \psi}$", which generalises other classical fractional derivatives, as follows
$$ \begin{aligned}
\mathbb{H}_p^{\alpha, \beta ; \psi} & =\left\{\begin{array}{c}
v \in L^p(\Lambda, \mathbb{R}) ;{ }^{\mathbf{H}} \mathbf{D}_{0+}^{\alpha, \beta ; \psi} v \in L^p(\Lambda, \mathbb{R}), \\
\mathbf{I}_{0+}^{\beta(\beta-1)} v(0)=\mathbf{I}_{T-}^{\beta(\beta-1)} v(T)=0
\end{array}\right\} & =\overline{\mathcal{C}_0^{\infty}(\Lambda, \mathbb{R})}
\end{aligned}
$$ 
equipped with the following norm
$$
\|v\|_{\mathbb{H}_p^{\alpha, \beta ; \psi}}=\left(\|v\|_{L^p}^p+\left\|^{\mathbf{H}} \mathbf{D}_{0+}^{\alpha, \beta ; \psi} v\right\|_{L^p}^p\right)^{1 / p},
$$
where $0<\alpha \leq 1$ and $0 \leq \beta \leq 1$.
For further details about the different parameters given above, see \cite{Costa}.\\
% Variable exponents in mathematics play a crucial role in a variety of fields, offering considerable power and flexibility for expressing and analysing complex mathematical concepts.
%In 2022 Sousa et al. have proposed a novel fractional space using $\psi$-Hilfer derivative incorporating the variable exponent, they defined the space $\mathcal{H}_{p(x)}^{\gamma, \beta ; \chi}(\Lambda)$ as the closure space of $\mathscr{C}_0^{\infty}(\Lambda)$ with the following norm
%$$
%\|v\|_{\mathcal{H}_{p(x)}^{\gamma, \beta ; \chi}}=\inf \left\{\lambda>0: \int_{\Lambda}\left|\frac{v}{\lambda}\right|^{p(x)}+\left|\frac{^{\mathrm{H}}{\mathcal{D}_{+}^{\gamma, \beta ; \chi}v}}{\lambda}\right|^{p(x)} d x \leqslant 1\right\},
%$$
%which is provided by
%$$
%\mathcal{H}_{p(x)}^{\gamma, \beta ; \chi}:=\left\{v \in L^{p(x)}(\Delta):{ }^{\mathrm{H}} \mathcal{D}_{+}^{\gamma, \beta ; \chi} v \in L^{p(x)}(\Lambda) \text { and } v(\Lambda)=0\right\} \text {. }
%$$

\noindent It is therefore a natural inquiry to explore the possibility of defining a new different type of fractional space using Orlicz spaces and fractional calculus. This exploration is valuable because extending classical results to new contexts can lead to a deeper and more comprehensive understanding of these results.\\  Thus, Inspired by the above interesting works,  our primary goal in this work is to define the new fractional Caputo-derivative Orlicz space  $\mathcal{O}^{\alpha,G}(\Lambda, \mathbb{R})$ and show several of there fundamental functional properties such as completeness,
reflexivity, separability, and establish some embeddings results.\\
Te remainder of this paper is abstracted as follows :
In Section \ref{sec2}, we recall the necessary theoretical background of Orlicz spaces, fractional calculus and some results that will be often used. In Section \ref{sec3}, we introduce a new fractional Caputo-derivative Orlicz space and demonstrate several of its properties. As an application, we use the mountain pass theorem to analyse a fractional boundary value problem.
\section{Preliminaries and basic notations}\label{sec2}
In this section, we sketch the basic facts of Orlicz spaces and fractional calculus.
\subsection{Orlicz spaces}
We devote the present subsection to listing briefly some fundamental terms and characteristics of Orlicz spaces. We would like to refer the reader to \cite{Adams,Rao,Mihailescu,Rutickii} and the references therein. \\
Let $\Omega$ be an open subset of $\mathbb{R}^N$. Assume that $G: \mathbb{R}^{+} \rightarrow \mathbb{R}^{+}$be a $\mathrm{N}$-function, meaning that $G$ is convex and continuous for all $s>0$, which satisfies $$
 \lim _{s \rightarrow 0^{+}} \frac{G(s)}{s}=0 \quad \text{ and } \quad
 \lim _{s \rightarrow \infty} \frac{G(s)}{s}=\infty.
$$ 
 Equivalently, $G$ admits the representation : $G(s)=\int_0^s g(t) d t$, where $g$ :
$\mathbb{R}^{+} \rightarrow \mathbb{R}^{+}$ satisfies the following conditions :
\begin{itemize}
    \item $ \lim _{s \rightarrow \infty} g(s)=\infty~$ and  $~g(0)=0 \text {; }$
    \item If $s>0$, then $g(s)>0$;
    \item $g$ is increasing;
    \item $g$ is right continuous.
\end{itemize}
We define the Orlicz space for the function G described before by
$$
L^{G}(\Omega)=\left\{v: \Omega \longrightarrow \mathbb{R} \text { measurable } \int_{\Omega} G(\gamma|v(s)|) d s<\infty, \text { for some } \gamma>0\right\}.
$$
The space $L^{G}(\Omega)$ is a Banach space equipped with the Luxemburg norm  
\begin{equation}\label{Luxemburg}
\|v\|_{G}=\inf \left\{\gamma>0: \int_{\Omega} G\left(\frac{|v(s)|}{\gamma}\right) d s \leqslant 1\right\} .
\end{equation}
The formula $\overline{G}(s)=\int_0^s \overline{g}(t) d t$ defines the conjugate $N$-function of $G$, where $\overline{g}: \mathbb{R}^+ \longrightarrow \mathbb{R}^+$ is given by $\overline {g}(s)=\sup \{t: g(t) \leqslant s\}$. \\ In addition, the Hölder inequality in Orlicz space is given by the following formula  
\begin{equation}\label{21}
    \left|\int_{\Omega} v u d s\right| \leqslant C\|v\|_{G}|| u \|_{\overline{G}} \quad \forall v \in L^{G}(\Omega) \text { and } u \in L^{\overline{G}}(\Omega),
\end{equation}
where $C$ is a positive constant (see {\cite[Theorem 13.13]{MusiElaK}}). \\ The function $G$ is said to satisfy the $\Delta_2$-condition (see \cite{Mihailescu})  if there exists $k>0$ such that  
\begin{equation}\label{delta}
G(2 s) \leqslant k G(s), \forall s \in \mathbb{R}_{+}.
\end{equation}

%When this inequality holds only for $t \geqslant t_0>0,$ $G$ is said to satisfy the $\Delta_2$-condition near infinity.
%We call the pair $(G, \Omega)$ is $\Delta$-regular if either :
%\begin{enumerate}
 %   \item[a)]  $G$ satisfies a global $\Delta_2$-condition, or
%\item[b)] $G$ satisfies a $\Delta_2$-condition near infinity and $\Omega$ has finite volume.
%
%\end{enumerate}
\noindent Given  an N-function $G$   which satisfies the global $\Delta_2$-condition. Then, we have
\begin{equation}\label{23:}
\bar{G}(g(s)) \leq c G(s) \text { for all } s \geq 0,
\end{equation}
where $c>0$.\\
Furthermore, the function $G$ has the following properties \cite{Mihailescu} :\\
%$\bullet$ (Subadditivity) Given $x, y \in \mathbb{R}^{+}$, $$
%G(x+y) \leq \frac{k}{2}(G(x)+G(y)),
%$$
%where $k>0$ is the constant in the $\Delta_2$-condition .\\
$\bullet$ For any $0<r<1$ and $s>0$, it holds 
\begin{equation}\label{23.}
    G(sr) \leqslant r G(s).
\end{equation}
$\bullet$ Every N-function $G$ satisfies
\begin{equation}\label{Gg}
    G(s) \leq s g(s) \leq G(2 s), \quad  \text{ for all } s \geq 0.
\end{equation}
In this paper we assume that
\begin{equation}\label{++}
    1<g^{-}:=\inf _{s>0} \frac{s g(s)}{G(s)}<g^{+}:=\sup _{s>0} \frac{s g(s)}{G(s)}<+\infty \text {. }
\end{equation}
In addition, we presume that $G$ satisfies the subsequent requirement the function
\begin{equation}\label{23}
    [0, \infty) \ni s \mapsto G(\sqrt{s}) \text{ is convex}.
\end{equation}
Relation (\ref{23}) assures that $L^{G}(\Omega)$ is an uniformly convex space consequently is a reflexive space (see {\cite[Proposition 2.2]{Mihailescu}}).\\
For all $ \phi \in L^{G}(\Lambda)$, we define the modular
$$
\rho(\phi)=\int_{\Lambda}G\left(|\phi(x)|\right) d x.
$$
If $\left(\phi_n\right), \phi \in L^{G}(\Lambda)$ then the relations shown below are valid
\begin{equation}\label{25}
    \|\phi\|_{G}>1 \Rightarrow\|\phi\|_{G}^{g^-} \leqslant \rho(\phi) \leqslant\|\phi\|_{G}^{g^+};
\end{equation}

\begin{equation}\label{26}
     \|\phi\|_{G}<1 \Rightarrow\|\phi\|_{G}^{g^+} \leqslant \rho(\phi) \leqslant\|\phi\|_{G}^{g^-};
\end{equation}

\begin{equation}
    \left\|\phi_n-\phi\right\|_{G} \rightarrow 0 \quad \Longleftrightarrow \quad \rho\left(\phi_n-\phi\right) \rightarrow 0;
\end{equation}

\begin{equation}
    \left\|\phi_n\right\|_{G} \rightarrow \infty \quad \Longleftrightarrow \quad \rho\left(\phi_n\right) \rightarrow \infty.
\end{equation}

\subsection{Fractional calculus} 
In the rest of this section, we define the Caputo fractional derivatives and the Riemann-Liouville fractional integrals.  For the proof see references \cite{Kilbas,Samko} and the reference  included.
\begin{definition}
    Define $(a,b)$ as a non-empty interval in $\mathbb{R}$ (with $-\infty \leqslant a<b \leqslant \infty$). The left and right-sided Riemann-Liouville fractional integrals $\mathbf{I}_{a+}^\alpha v$ and $\mathbf{I}_{b-}^\alpha v$ of order $\alpha \in \mathbb{R}_{+}$ are defined by
     \begin{equation}
\mathbf{I}_{a+}^\alpha v(x):=\frac{1}{\Gamma(\alpha)} \int_a^x v(t)(x-t)^{\alpha-1} d t
 \end{equation}
and
\begin{equation}
\mathbf{I}_{b-}^\alpha v(x):=\frac{1}{\Gamma(\alpha)} \int_x^b v(t)(t-x)^{\alpha-1}  d t
 \end{equation}
 respectively. Where $\Gamma(\alpha)$ is the Gamma function. 
 \end{definition}
 \begin{definition}
 The Riemann-Liouville fractional derivatives $D_{a+}^\alpha v$ and $D_{b-}^\alpha v$ of order $\alpha \in \mathbb{R}_+$ are defined by
\begin{equation}\label{210}
    \begin{aligned}
D_{a+}^\alpha v(x):&=\left(\frac{d}{d x}\right)^n\mathbf{I}_{a+}^{n-\alpha} v(x) \\
&=\frac{1}{\Gamma(n-\alpha)}\left(\frac{d}{d x}\right)^n \int_a^x v(t)(x-t)^{n-\alpha-1} dt  
\end{aligned}
\end{equation}
and 
\begin{equation}\label{211}
    \begin{aligned}
D_{b-}^\alpha v(x)&:=\left(-\frac{d}{d x}\right)^n\mathbf{I}_{b-}^{n-\alpha} v(x) \\
&=\frac{1}{\Gamma(n-\alpha)}\left(-\frac{d}{d x}\right)^n \int_x^b v(t)(t-x)^{n-\alpha-1} dt,
\end{aligned}
\end{equation}
respectively, where $n=[\alpha]+1$.
 \end{definition}
\begin{definition}
 The left-sided and right-sided Caputo fractional derivatives ${ }^C \mathfrak{D}_{a+}^\alpha (\cdot)$ and ${ }^C \mathfrak{D}_{b-}^\alpha (\cdot)$ of order $\alpha $ on $[a, b]$ are defined via the above Riemann-Liouville fractional derivatives by
\begin{equation}\label{212}
{ }^C \mathfrak{D}_{a+}^\alpha v(x):=D_{a+}^\alpha\left[v(t)-\sum_{k=0}^{n-1} \frac{v^{(k)}(a)}{k !}(t-a)^k\right](x)
\end{equation}
and
\begin{equation}\label{213}
{ }^C \mathfrak{D}_{b-}^\alpha v(x):=D_{b-}^\alpha\left[v(t)-\sum_{k=0}^{n-1} \frac{v^{(k)}(b)}{k !}(b-t)^k\right](x)
\end{equation}
respectively, where  $D_{a+}^\alpha (\cdot)$ and $D_{b-}^\alpha (\cdot)$ are defined in $(\ref{210})$ and  $(\ref{211})$.
%and 
%\begin{equation}\label{214}
%n=[\alpha]+1 \text{ for } \alpha \notin \mathbb{N} ; \quad n=\alpha \text{ for } \alpha \in \mathbb{N}.
%\end{equation}
In particular, when $0<\alpha<1$, the relations $(\ref{212})$ and $(\ref{213})$ take the following forms:
$$
\begin{aligned}
& { }^C \mathfrak{D}_{a+}^\alpha v(x)=D_{a+}^\alpha[v(t)-v(a)](x), \\
& { }^C \mathfrak{D}_{b-}^\alpha v(x)=D_{b-}^\alpha[v(t)-v(b)](x) .
\end{aligned}
$$
\end{definition}
\begin{theorem}\label{th2.1}
Let $\alpha \geq 0$.
%and let $n$ be given by $(\ref{214})$. 
If $v\in A \mathscr{C}^n([a, b])$, then the Caputo fractional derivatives ${ }^C \mathfrak{D}_{a+}^\alpha (\cdot)$ and ${ }^C \mathfrak{D}_{b-}^\alpha ( \cdot)$ exist almost everywhere on $[a, b]$. Moreover
\begin{equation}
    { }^C \mathfrak{D}_{a+}^\alpha v(x)=\frac{1}{\Gamma(n-\alpha)} \int_a^x v^{(n)}(t) (x-t)^{n-\alpha-1} d t=:\mathbf{I}_{a+}^{n-\alpha} D^n v(x)
\end{equation}
and 
\begin{equation}\label{2.20}
    { }^C \mathfrak{D}_{b-}^\alpha v(x)=\frac{(-1)^n}{\Gamma(n-\alpha)} \int_x^b v^{(n)}(t) (t-x)^{n-\alpha-1} d t=:(-1)^n\mathbf{I}_{b-}^{n-\alpha} D^n v(x),
\end{equation}
where $D=d / d x$,  $n=\alpha$   for  $\alpha \in \mathbb{N}$ and  $n=[\alpha]+1$  for  $\alpha \notin \mathbb{N}$.
\end{theorem}
\noindent In particular, when $\alpha\in (0,1)$ and $v \in A \mathscr{C}([a, b])$, we have
\begin{equation}
    { }^C \mathfrak{D}_{a+}^\alpha v(x)=\frac{1}{\Gamma(1-\alpha)} \int_a^x v^{\prime}(t) (x-t)^{-\alpha} d t=:\mathbf{I}_{a+}^{1-\alpha} D v(x)
\end{equation}
and 
\begin{equation}
    { }^C \mathfrak{D}_{b-}^\alpha v(x)=-\frac{1}{\Gamma(1-\alpha)} \int_x^b v^{\prime}(t) (t-x)^{-\alpha} d t=:-\mathbf{I}_{b-}^{1-\alpha} D v(x).
\end{equation}
\begin{proposition}\label{P21}
     Let $\alpha>0$ and given $n$ as in \eqref{2.20}. If $v \in A \mathscr{C}^n([a, b])$ or $v \in \mathscr{C}^n\left( [a, b]\right) $, then
     \begin{equation}
        \mathbf{I}_{a+}^\alpha{ }^C \mathfrak{D}_{a+}^\alpha v(x)=v(x)-\sum_{k=0}^{n-1} \frac{v^{(k)}(a)}{k !}(x-a)^k
     \end{equation}
     and
     \begin{equation}
         \mathbf{I}_{a+}^\alpha{ }^C \mathfrak{D}_{a+}^\alpha v(x)=v(x)-\sum_{k=0}^{n-1} \frac{(-1)^k v^{(k)}(b)}{k !}(b-x)^k .
     \end{equation}
     In particular, if $\alpha \in (0,1)$ and $v \in A  \mathscr{C}([a, b])$ or $v \in \mathscr{C}([a, b])$, then
     $$\mathbf{I}_{a+}^\alpha{ }^C \mathfrak{D}_{a+}^\alpha v(x)=v(x)-v(a), \text{ and }  \mathbf{I}_{a+}^\alpha{ }^C \mathfrak{D}_{a+}^\alpha v(x)=v(x)-v(b).$$
\end{proposition}
\section{Fractional Caputo-derivative Orlicz  space}\label{sec3}
The present section is devoted to the new fractional Caputo-derivative Orlicz functional framework and the establishment of there basic properties. \\

 \noindent Given $0<\alpha< 1$ and let $G$ be an $N$-function and we denote $\Lambda=[0, T]$. Then, we introduce the left-sided fractional Caputo-derivative Orlicz  space $\mathcal{O}^{\alpha,G}(\Lambda, \mathbb{R})=\mathcal{O}^{\alpha,G}$ as follows :
    \begin{equation}\label{O}
        \mathcal{O}^{\alpha,G}=\mathcal{O}=\left\{v\in L^{G}(\Lambda, \mathbb{R}) ~~\text{ such that }~~ { }^C \mathfrak{D}_{0^+}^\alpha v\in L^{G}(\Lambda, \mathbb{R}) \right\}  .
        \end{equation}
We equip this space  with the following  norm 
\begin{equation}\label{32}
    \|v\|_{\mathcal{O}}=\|v\|_G+[v]_{\mathcal{O}},
\end{equation}
where $[\cdot]_{\mathcal{O}}$ is  defined by
\begin{equation}
    [v]_{\mathcal{O}}=\inf\left\{\gamma>0: \int_\Lambda G\left(\frac{\left|  { }^C \mathfrak{D}_{0^+}^\alpha v(x)\right|}{\gamma}\right) d x \leqslant 1  \right\}\text{.}
\end{equation}
\begin{remark}
$\mathcal{O}$ is non-empty. Indeed, let us consider the constant function $v(x)=C$,  $G(x)=e^x$ and let $\alpha=\tfrac{1}{2}$.\\
Clearly $v \in L^G(\Lambda)$, it remains to show that $ {}^C \mathfrak{D}_{0^+}^\alpha v\in L^{G}(\Lambda).$
 $$ {}^C \mathfrak{D}_{0^+}^\alpha v(x)=\frac{1}{\Gamma(\frac{1}{2})}\int_0^T (x-t)^{-\alpha} v^{\prime}(t)dt =0.$$
 Then $$\psi(x)={}^C \mathfrak{D}_{0^+}^\alpha v(x)=0.$$
Hence, we have $\psi \in L^G(\Lambda)$, which implies that $\mathcal{O}$ is non-empty.\\
Note that The fractional derivative of a constant is not necessarily zero. However, the Caputo derivative of a constant is always zero, which is a distinctive feature compared to other types of fractional derivatives.  
\end{remark}

\noindent Next, for any $v \in \mathcal{O}$, we define the modular function 
$$
\rho_{G}(v)=\int_{\Lambda} \left[G(|v(x)|)+G\left(\left|  { }^C \mathfrak{D}_{0^+}^\alpha v(x)\right|\right) \right]d x.
$$
It is obvious to check that $\rho_G$ is a convex modular on $\mathcal{O}$. Moreover, the  norm $(\ref{32})$ can be defined as follows:
\begin{equation}\label{33} 
\|v\|=\inf\left\{\gamma>0:\rho_{G}\left(\frac{v}{\gamma}\right) \leqslant 1\right\}.
\end{equation}
We also define the closure of $\mathscr{C}_0^{\infty}\left(\Lambda, \mathbb{R}\right)$, in the norm $\|\cdot\|$ defined in (\ref{33}) by 
$$\mathcal{O}_0=\left\{v\in \mathcal{O}: v(0)=v(T)=0\right\}.$$
On the other hand, for any $v\in \mathcal{O}_0$, we define the convex modular as follows
$$
\rho_G^0(v)=\int_{\Lambda} G\left(\left|  { }^C \mathfrak{D}_{0^+}^\alpha v(x)\right|\right) d x.
$$
The norm associated with $\rho_G^0$ is given by 
$$\|v\|_{0}=[v]_{\mathcal{O}}=\inf\left\{\gamma>0:\rho_G^0\left(\frac{v}{ \gamma}\right)\leqslant 1\right\}.$$
 \noindent Note that in the case $g(t)=p|t|^{p-2} t$, the space $\mathcal{O}_{0}$ coincide with the fractional derivative space $E_0^{\alpha, p}$ introduced in {\cite[Definition 6.1]{Zhou}}.

\begin{theorem}
    The space $(\mathcal{O},\|\cdot\|_\mathcal{O})$ defined by $(\ref{O})$ is a Banach space.
\end{theorem}
\noindent \textbf{\textit{Proof.}} 
We define the operator $\mathcal{T}$ by
$$
\begin{aligned}
\mathcal{T}: \mathcal{O} & \longrightarrow L^{G}(\Lambda) \times  L^{G}(\Lambda)=\Pi \\
v & \longmapsto\left(v(x), { }^C \mathfrak{D}_{0^+}^\alpha v(x)\right) .
\end{aligned}
$$
For any $v \in \mathcal{O}$, we have
$$\|\mathcal{T}(v)\|_{\Pi}=\|v\|_G+\|{ }^C \mathfrak{D}_{0^+}^\alpha v\|_G=\|v\|_G+[v]_{\mathcal{O}}=\|v\|_{\mathcal{O}}.$$
Then $\mathcal{T}$ is an isometry. Since $\Pi$ is a Banach space \cite{Adams}, therefore $\mathcal{O}$ is also a Banach space.
\hfill$\Box$
\begin{lemma}
    Let $0<\alpha <1$ and $G$ be an $N$-function, we suppose that conditions $(\ref{delta})$ and $(\ref{23})$ are satisfied. Then $(\mathcal{O},\|\cdot\|_{\mathcal{O}})$ is separable and uniformly convex space (consequently it is a reflexive space). 
\end{lemma}
\noindent \textbf{\textit{Proof.}}  Using the operator $\mathcal{T}$ introduced in the previous theorem, we get that $\mathcal{T}(\mathcal{O})$ is closed subset of $\Pi$, since $\mathcal{O}$ is a Banach space. Consequently, by ({\cite[Proposition 3.25]{Brezis}}), the space $(\mathcal{O},\|\cdot\|_{\mathcal{O}})$ is separable.\\
On the other hand, by ({\cite[Theorem 11.6]{MusiElaK}}) the space $(L^G,\|\cdot\|_G)$ is uniformly convex. Moreover by ({\cite[page 184]{MusiElaK}}), we deduce that $\Pi$ is uniformly convex. Since $\mathcal{T}(\mathcal{O})$ is a closed subset of $\Pi$ and any linear subspace of a uniformly convex linear normed space is also uniformly convex. Hence $(\mathcal{O},\|\cdot\|_{\mathcal{O}})$ is uniformly convex, as a result, it is a reflexive space. The proof is completed.$\hspace{0.5 cm}\Box$
\begin{proposition}\label{P32}
    Let $v\in \mathcal{O}$ and assume that $(\ref{++})$ is satisfied. Then the following relations hold true
\begin{equation}\label{35}
\bullet~~  \text{ If } [v]_{\mathcal{O}}>1, \text{ then } [v]_{\mathcal{O}}^{g^{-}} \leqslant \rho_G^0(v) \leqslant[v]_{\mathcal{O}}^{g^{+}}, \text{ for all } v \in \mathcal{O}  .
    \end{equation}
   \begin{equation}\label{36}
\bullet~~   \text{ If } [v]_{\mathcal{O}}<1, \text{ then } [v]_{\mathcal{O}}^{g^{+}} \leqslant \rho_G^0(v) \leqslant[v]_{\mathcal{O}}^{g^{-}}, \text{ for all } v \in \mathcal{O} .
    \end{equation}
    
   % \begin{equation}\label{37}
      %\left[u_n-u\right]_\mathcal{O} \rightarrow 0 \quad \Leftrightarrow \quad \rho_G^0\left(u_n-u\right) \rightarrow 0
  %\end{equation}

  %\begin{equation}\label{38}
      %\left[u_n\right]_{\mathcal{O}} \rightarrow \infty \quad \Leftrightarrow \quad \rho_G^0\left(u_n\right) \rightarrow \infty
  %\end{equation}
\end{proposition}
\noindent \textbf{\textit{Proof.}} We follow her an approach similar to the one used in \cite{Mihailescu}. 
Initially, we demonstrate that  $\rho_G^0(v) \leqslant[v]_{\mathcal{O}}^{g^{+}},$  for all  $v \in \mathcal{O} \text{ with } [v]_{\mathcal{O}}>1 .$ Indeed, since $g^{+} \geqslant \frac{s g(s)}{G(s)}$ for all $s>0$, it follows that for letting $r>1$ we obtain 
$$
\log (G(r s))-\log (G(s))=\int_s^{r s} \frac{g(\tau)}{G(\tau)} d \tau \leqslant \int_s^{r s} \frac{g^{+}}{\tau} d \tau=\log \left(r^{g^{+}}\right) .
$$
Hence, we infer  
\begin{equation}\label{37}
    G(r s) \leqslant r^{g^{+}} G(s) \text { for all } s>0 \text { and } r>1 \text {. }
\end{equation}
Let now $v \in \mathcal{O}$ with $[v]_{\mathcal{O}}>1$. Using the relation $(\ref{37})$ and the definition of the Luxemburg norm, we get
$$
\begin{aligned}
\int_{\Lambda}  G\left(\left|{ }^C \mathfrak{D}_{0^+}^\alpha v(s)\right|\right) d s & = \int_{\Lambda} G\left([v]_\mathcal{O} \tfrac{\left|{ }^C \mathfrak{D}_{0^+}^\alpha v(s)\right|}{[v]_\mathcal{O}}\right) d s \\
& \leqslant[v]_\mathcal{O}^{g^{+}} \int_{\Lambda}  G\left(\tfrac{\left|{ }^C \mathfrak{D}_{0^+}^\alpha v(s)\right|}{[v]_{\mathcal{O}}}\right) d s \\
& \leqslant[v]_{\mathcal{O}}^{g^{+}} .
\end{aligned}
$$
Now, we show that  $ [v]_{\mathcal{O}}^{g^{-}} \leqslant \rho_G^0(v) $ for all $v \in \mathcal{O}$ with $[v]_{\mathcal{O}}>1$. By applying a method similar to the one previously established in equation $(\ref{37})$, we deduce
\begin{equation}\label{38}
    G(r s) \geqslant r^{g^{-}} G(s) \text { for all } s>0 \text { and } r>1 \text {. }
\end{equation}
Let $v \in\mathcal{O}$ with $[v]_{\mathcal{O}}>1$, we consider $\sigma \in\left(1,[v]_{\mathcal{O}}\right)$, since $\sigma<[v]_{\mathcal{O}}$. Thus, per the Luxemburg norm's definition, it follows that
 $ \int_{\Lambda} G\left(\frac{\left|{ }^C \mathfrak{D}_{0^+}^\alpha v(s)\right|}{\sigma}\right) d s >1 $, if not we will obtain a contradiction with the definition of the Luxemburg norm. The idea mentioned above suggests that 
$$
\begin{aligned}
\int_{\Lambda} G\left(\left|{ }^C \mathfrak{D}_{0^+}^\alpha v(s)\right|\right) d s & =\int_{\Lambda}  G\left(\sigma \frac{\left|\mathfrak{D}_{0^+}^\alpha v(s)\right|}{\sigma}\right) d s \\
& \geqslant \sigma^{g^{-}}  \int_{\Lambda} G\left(\frac{\left|\mathfrak{D}_{0^+}^\alpha v(s)\right|}{\sigma}\right) d s \\
& \geqslant \sigma^{g^{-}},
\end{aligned}
$$
 when $\sigma$ tends to $[v]_{\mathcal{O}}$, we conclude that relation $(\ref{35})$ is valid.\\
Next, we prove that  $\rho_G^0(v) \leqslant[v]_{\mathcal{O}}^{g^{-}} \text{ for all } v \in \mathcal{O} \text{ with } [v]_{\mathcal{O}}<1.$ Utilising a method similar to the one previously given in relation $(\ref{37})$, we obtain 
\begin{equation}\label{39}
    G(s) \leqslant \tau^{g^{-}} G\left(\tfrac{s}{\tau}\right),~~ \text { for all } s>0 \text { and } \tau \in(0,1) \text {. }
\end{equation}
Let $v \in \mathcal{O}$ with $[v]_{\mathcal{O}}<1$. From the definition of the norm define by $(\ref{Luxemburg})$ and the relation $(\ref{39})$, we conclude  
$$
\begin{aligned}
 \int_{\Lambda} G\left(\left|{ }^C \mathfrak{D}_{0^+}^\alpha v(s)\right|\right) d s & =\int_{\Lambda} G\left(\frac{\left|{ }^C \mathfrak{D}_{0^+}^\alpha v(s)\right| [v]_{\mathcal{O}} }{[v]_{\mathcal{O}}}\right) d s\\
 &\leqslant[v]_{\mathcal{O}}^{g^{-}}  \int_{\Lambda} G\left(\frac{\left|{ }^C \mathfrak{D}_{0^+}^\alpha v(s)\right|}{[v]_{\mathcal{O}}}\right) d s \\
& \leqslant[v]_{\mathcal{O}}^{g^{-}} .
\end{aligned}
$$
Finally, we establish that $[v]_{\mathcal{O}}^{g^{+}} \leqslant \rho_G^0(v) \text{ for all } v \in \mathcal{O} \text{ with } [v]_{\mathcal{O}}<1.$ Similar manner as the one used in the proof of relation $(\ref{37})$, we have 
\begin{equation}\label{310}
    G(s) \geqslant \tau^{g^{+}} G\left(\frac{s}{\tau}\right) \text { for all } s>0 \text { and } \tau \in(0,1) \text {. }
\end{equation}
Let $v \in \mathcal{O}$ with $[v]_{\mathcal{O}}<1$ and $\sigma \in\left(0,[v]_{\mathcal{O}}\right)$, so by $(\ref{310})$ we find
\begin{equation}\label{311}
    \int_{\Lambda} G\left(\left|{ }^C \mathfrak{D}_{0^+}^\alpha v(s)\right|\right) d s \geqslant \sigma^{g^{+}} \int_{\Lambda}  G\left(\frac{\left|{ }^C \mathfrak{D}_{0^+}^\alpha v(s)\right|}{\sigma}\right) d s .
\end{equation}
We define $u(x)=\frac{v(s)}{\sigma}$ for all $s \in \Lambda$, we have $[u]_{\mathcal{O}}=\frac{[v]_{\mathcal{O}}}{\sigma}>1$.
Using the relation $(\ref{35})$, we find 
\begin{equation}\label{312}
 \int_{\Lambda} G\left(\frac{\left|{ }^C \mathfrak{D}_{0^+}^\alpha v(s)\right|}{\sigma}\right) d s=\int_{\Lambda}  G\left(\left|{ }^C \mathfrak{D}_{0^+}^\alpha u(s)\right|\right) d s>[u]_{\mathcal{O}}^{g^{-}}>1.
\end{equation}
By $(\ref{311})$ and $(\ref{312})$, we obtain
$$
 \int_{\Lambda} G\left(\left|{ }^C \mathfrak{D}_{0^+}^\alpha v(x)\right|\right) d x \geqslant \sigma^{g^{+}} .
$$
Letting $\sigma \nearrow[v]_{\mathcal{O}}$, we deduce that relation $(\ref{36})$ hold true.
  %One can infer relationships $(\ref{35})$ and $(\ref{36})$ directly from relationships $(\ref{37})$ and $(\ref{38})$.  
\hfill$\Box$\\

\noindent In the outcomes that follow, for $\alpha\in(0,1)$, we presume that 
\begin{equation}\label{313.}
    (x -\tau)^{\alpha}<(x -\tau), 
 \hspace{0.3 cm}\text{ for all } x\in \Lambda\text{ and } \tau\in [0,x ]. 
\end{equation}
and we utilise $g^{\pm}$ instead of recalculating, where $g^+$ or $g^-$ is selected based on whether $[u]$ is higher   or lower than 1. 
\begin{proposition}\label{P33}
    Let $0<\alpha<1$ and $G$ be an N-function. Assume that $(\ref{313.})$ be verified, then for any $v\in L^G(\Lambda)$, we have 
    $$
\left\{\begin{array}{l}
\left\|\mathbf{I}_{0^+}^{\alpha} v\right\|_{G} \leqslant \left[\tfrac{T^{\alpha}}{\Gamma(\alpha+1)}\right]^{\frac{1}{g^-}} \|v\|_{G}^{\frac{g^+}{g^-}},~~ \text{ if } \|\cdot\|>1,\\
\left\|\mathbf{I}_{0^+}^{\alpha} v\right\|_{G} \leqslant \left[\tfrac{T^{\alpha}}{\Gamma(\alpha+1)}\right]^{\frac{1}{g^+}} \|v\|_{G}^{\frac{g^-}{g^+}},~~ \text{ if } \|\cdot\|<1,\\
\end{array}\right.
$$
where $\|\cdot\|>1$ means that $\|v\|_G>1$ and /or $\|\mathbf{I}_{0^+}^{\alpha} v\|_G>1$, we use the same notation for to the lower case $(<)$.
\end{proposition}
\noindent \textbf{\textit{Proof.}} 
By Dirichlet formula and Jensen's inequality, we obtain 
$$\begin{aligned}
\rho_G^0\left(\mathbf{I}_{0^+}^{\alpha} v(x)\right)
& =\int_0^T G\left(\left|\frac{1}{\Gamma(\alpha)} \int_0^{x} v(\tau)(x-\tau)^{\alpha-1} d \tau\right|\right) dx \\
& \leqslant \int_0^T \int_0^{x} G\left(\left|\frac{1}{\Gamma(\alpha)}  v(\tau)(x-\tau)^{\alpha-1} d \tau\right|\right) dx\\
& \leqslant \frac{1}{\Gamma(\alpha)} \int_0^T \int_0^{x} (x-\tau)^{\alpha-1} G\left(\left|  v(\tau) \right|\right) d x d \tau  \\
&= \frac{1}{\Gamma(\alpha)} \int_0^T G\left(\left|  v(\tau)\right|\right) \int_{\tau}^{T} (x-\tau)^{\alpha-1}   dx d \tau =\frac{T^{\alpha}}{\Gamma(\alpha+1)} \rho_G^0(v)
\end{aligned}$$
Hence, if $\|\cdot\|>1$ by relation $(\ref{25})$, we have 
$$\left\|\mathbf{I}_{0^+}^{\alpha} v\right\|_{G} \leqslant \left[\tfrac{T^{\alpha}}{\Gamma(\alpha+1)}\right]^{\frac{1}{g^-}} \|v\|_{G}^{\frac{g^+}{g^-}}.$$
Similarly, if $\|\cdot\|<1$ by relation $(\ref{26})$, we have 
$$\left\|\mathbf{I}_{0^+}^{\alpha} v\right\|_{G} \leqslant \left[\tfrac{T^{\alpha}}{\Gamma(\alpha+1)}\right]^{\frac{1}{g^+}} \|v\|_{G}^{\frac{g^-}{g^+}}.$$
\hfill$\Box$
\begin{proposition}\label{P34}
    Let $0<\alpha<1$ and $G$ be an N-function. Assume that $(\ref{313.})$ de satisfied. For all $v\in \mathcal{O}_0$, we have $ \mathbf{I}_{0^+}^\alpha{ }^C \mathfrak{D}_{0^+}^\alpha v(t)=v(t)$. Moreover the inclusion $\mathcal{O}_0\subset \mathscr{C}(\Lambda)$ holds. 
\end{proposition}
\noindent \textbf{\textit{Proof.}}  For any $0<t_1<t_2\leqslant T$, using $(\ref{21})$ and $(\ref{23.})$, we have 
$$
{\footnotesize \begin{aligned}
\left|\mathbf{I}_{0^+}^\alpha v\left(t_2\right)-\mathbf{I}_{0^+}^\alpha v\left(t_1\right)\right| & =\frac{1}{\Gamma(\alpha)}\left|\int_0^{t_2}\left(t_2-\tau\right)^{\alpha-1} v(\tau) d \tau-\int_0^{t_1}\left(t_1-\tau\right)^{\alpha-1} v(\tau) d \tau\right| \\
& \leqslant \frac{1}{\Gamma(\alpha)}\left|\int_{t_1}^{t_2}\left(t_2-\tau\right)^{\alpha-1} v(\tau) d \tau\right| \\
& +\frac{1}{\Gamma(\alpha)}\left|\int_0^{t_1}\left(\left(t_2-\tau\right)^{\alpha-1}-\left(t_1-\tau\right)^{\alpha-1}\right) v(\tau) d \tau\right|\\
& \leqslant \frac{C}{\Gamma(\alpha)} \left\|v \right\|_{G} \left\|\left(t_2-\tau\right)^{\alpha-1}\right\|_{\overline{G}}+ \frac{C}{\Gamma(\alpha)} \left\|v \right\|_{G} \left\|\left(\left(t_2-\tau\right)^{\alpha-1}-\left(t_1-\tau\right)^{\alpha-1}\right) \right\|_{\overline{G}} \\
& \leqslant \frac{C}{\Gamma(\alpha)} \left\|v \right\|_{G}  \left(\int_{t_1}^{t_2} \overline{G}\left(\left(t_2-\tau\right)^{\alpha-1}\right)d\tau\right)^{\frac{1}{g^+}}\\
& +\frac{C}{\Gamma(\alpha)} \left\|v \right\|_{G} \left(\int_{t_1}^{t_2} \overline{G}\left(\left(t_2-\tau\right)^{\alpha-1}-\left(t_1-\tau\right)^{\alpha-1}\right)d\tau\right)^{\frac{1}{g^+}} \\
&\leqslant \frac{C}{\Gamma(\alpha)} \left\|v \right\|_{G} \times  \left(\int_{t_1}^{t_2} \overline{G}(1)\left(t_2-\tau\right)^{\alpha-1}d\tau\right)^{\frac{1}{g^+}}  \\ 
 &+ \frac{C}{\Gamma(\alpha)} \left\|v \right\|_{G} \left(\int_{0}^{t_1} \overline{G}(1)\left(\left(t_2-\tau\right)^{\alpha-1}-\left(t_1-\tau\right)^{\alpha-1}\right)d\tau\right)^{\frac{1}{g^+}} 
 \\
& \leqslant \frac{C\left(\overline{G}(1)\right)^\frac{1}{g^+}}{\Gamma(\alpha)} \left\|v \right\|_{G} \\
&\times \left( \left(\int_{t_1}^{t_2} \left(t_2-\tau\right)^{\alpha-1}d\tau\right)^{\frac{1}{g^+}}    
 +\left(\int_{0}^{t_1} \left(\left(t_2-\tau\right)^{\alpha-1}-\left(t_1-\tau\right)^{\alpha-1}\right)d\tau\right)^{\frac{1}{g^+}} \right)
\end{aligned}}$$
\begin{equation}\label{313}
{\footnotesize \begin{aligned}
\left|\mathbf{I}_{0^+}^\alpha v\left(t_2\right)-\mathbf{I}_{0^+}^\alpha v\left(t_1\right)\right| &\leq \frac{C\left(\overline{G}(1)\right)^\frac{1}{g^+}}{\Gamma(\alpha+1)}\left\|v \right\|_{G} \left(\left(t_2-t_1\right)^{\frac{\alpha}{g^+}}+\left[ (t_2-0)^{\alpha}-(t_2-t_1)^{\alpha}-(t_1-0)^{\alpha} \right]^{\frac{1}{g^+}}\right)\\
 &\leqslant \frac{3C\left(\overline{G}(1)\right)^\frac{1}{g^+}}{\Gamma(\alpha+1)}\left\|v \right\|_{G}(t_2-t_1)^{\frac{\alpha}{g^+}}.\end{aligned}}
\end{equation}
Therefore 
$$\left|\mathbf{I}_{0^+}^\alpha v\left(t_2\right)-\mathbf{I}_{0^+}^\alpha v\left(t_1\right)\right| \leqslant \frac{3C\left(\overline{G}(1)\right)^\frac{1}{g^+}}{\Gamma(\alpha+1)}\left\|v \right\|_{G}(t_2-t_1)^{\frac{\alpha}{g^+}}. $$
From Proposition $\ref{P21}$, we have 
$$\mathbf{I}_{0^+}^\alpha{ }^C \mathfrak{D}_{0^+}^\alpha v(t)=v(t)-v(0).$$
Since ${ }^C \mathfrak{D}_{0^+}^\alpha v\in L^G$ by $(\ref{313})$, we obtain the continuity of $\mathbf{I}_{0^+}^\alpha{ }^C \mathfrak{D}_{0^+}^\alpha v$ in $\Lambda$ and since $v(0)=0$, then 
$ \mathbf{I}_{0+}^\alpha{ }^C \mathfrak{D}_{0+}^\alpha v(t)=v(t)$.  \hfill$\Box$
\begin{proposition}\label{P35}
    Let $0<\alpha<1$ and $G$ be an N-function. Assume that $(\ref{313.})$ verified. Then for all $v\in \mathcal{O}_0$, we have 
    \begin{equation}\label{315.}
\left\{\begin{array}{l}
\left\| v\right\|_{G} \leqslant \left[\frac{T^{\alpha}}{\Gamma(\alpha+1)}\right]^{\frac{1}{g^-}} \left[v\right]_{\mathcal{O}}^{\frac{g^+}{g^-}},~~ \text{ if } \|\cdot\|>1,\\
\left\| v\right\|_{G} \leqslant \left[\frac{T^{\alpha}}{\Gamma(\alpha+1)}\right]^{\frac{1}{g^+}} \left[ v\right]_{\mathcal{O}}^{\frac{g^-}{g^+}},~~ \text{ if } \|\cdot\|<1.\\
\end{array}\right.
\end{equation}   
    Moreover 
    \begin{equation}\label{IN}
        \|v\|_{\infty}\leqslant \frac{C\left(\overline{G}(1)\right)^\frac{1}{g+}}{\Gamma(\alpha+1)}x^{\alpha} \left[ v\right]_{\mathcal{O}}.
        \end{equation}
\end{proposition}
\noindent \textbf{\textit{Proof.}}  Since ${ }^C \mathfrak{D}_{0^+}^\alpha v \in L^G$, it follows from Proposition $\ref{P33}$ that 
   $$
\left\{\begin{array}{l}
\left\| \mathbf{I}_{0^+}^\alpha{ }^C \mathfrak{D}_{0^+}^\alpha v\right\|_{G} \leqslant \left[\tfrac{T ^{\alpha}}{\Gamma(\alpha+1)}\right]^{\frac{1}{g^-}} \left[v\right]_{\mathcal{O}}^{\frac{g^+}{g^-}},~~ \text{ if } \|\cdot\|>1,\\
\left\| \mathbf{I}_{0^+}^\alpha{ }^C \mathfrak{D}_{0^+}^\alpha v\right\|_{G}  \leqslant \left[\tfrac{T^{\alpha}}{\Gamma(\alpha+1)}\right]^{\frac{1}{g^+}} \left[ v\right]_{\mathcal{O}}^{\frac{g^-}{g^+}},~~ \text{ if } \|\cdot\|<1.\\
\end{array}\right.
$$   
Using the Proposition $\ref{P34}$, we obtain the first result $(\ref{315.})$. \\
By $(\ref{21})$ and $(\ref{23.})$, we have for all $v\in \mathcal{O}_0$ that 
$$
\begin{aligned}
    \left| \mathbf{I}_{0^+}^\alpha{ }^C \mathfrak{D}_{0^+}^\alpha v\right|&\leqslant\left|\frac{1}{\Gamma(\alpha)} \int_0^x (x-t)^{\alpha-1} { }^C \mathfrak{D}_{0^+}^\alpha v(t) dt \right|\\
    & \leqslant \frac{1}{\Gamma(\alpha)} \int_0^x (x-t)^{\alpha-1} \left|{ }^C \mathfrak{D}_{0^+}^\alpha v(t)\right| dt\\
    &\leqslant \frac{C}{\Gamma(\alpha)} \left[ v\right]_{\mathcal{O}}  \left\|(x-t)^{\alpha-1}\right\|_{\overline{G}}\\
    &\leqslant \frac{C}{\Gamma(\alpha)} \left[ v\right]_{\mathcal{O}} \left(\int_0^x\overline{G}\left((x-t)^{\alpha-1}\right) dt \right)^{\frac{1}{g^+}}\\
    &\leqslant \frac{C}{\Gamma(\alpha)} \left[ v\right]_{\mathcal{O}} \left(\int_0^x\overline{G}(1)(x-t)^{\alpha-1} dt \right)^{\frac{1}{g^+}} = \frac{C\left(\overline{G}(1)\right)^{\frac{1}{g^+}}}{\Gamma(\alpha+1)} x^{\frac{\alpha}{g}}\left[ v\right]_{\mathcal{O}} .
\end{aligned}
$$\hfill$\Box$

Note that from $(\ref{315.})$, we can to take in consideration on $\mathcal{O}_0$ the following norm 
\begin{equation}\label{317.}
\|v\|_{\mathcal{O}_0}=
\|{}^C \mathfrak{D}_{0^+}^\alpha v\|_G.
\end{equation}Moreover, $(\ref{317.})$ is equivalent to $(\ref{33})$ on $\mathcal{O}_0$.

\begin{proposition}
    Let $0<\alpha<1$ . Assume that $(\ref{313.})$ be satisfied. Then the embedding $\mathcal{O}_0 \hookrightarrow \mathscr{C}(\Lambda)$ is compact.
\end{proposition}
\noindent \textbf{\textit{Proof.}} Since $\mathcal{O}_0$ is a reflexive Banach space, it suffices to show that for any $\left(\varphi_n\right)_{n \in \mathbb{N}}$ converges weakly  to $\varphi$ in $\mathcal{O}_0$, then $\varphi_n$ converges strongly in $\mathscr{C}$.\\
Indeed, given $\left(\varphi_n\right)_{n \in \mathbb{N}} \subset \mathcal{O}_0$ such that $$\varphi_n \stackrel{\mathcal{O}_0}{\rightharpoonup} \varphi.$$ As $\mathcal{O}_0 \hookrightarrow \mathscr{C}$, we have $$\varphi_n \stackrel{\mathscr{C}}{\rightharpoonup} \varphi.$$
Hence $\left(\varphi_n\right)_{n \in \mathbb{N}}$ is bounded in $\mathcal{O}_0$ since it converges weakly in that space. Thus, $\left(^C{\mathfrak{D}_{0^+}^{\alpha}} \varphi_n\right)_{n \in \mathbb{N}}$ is bounded in $L^{G}$. In the following, we show that $\left(\varphi_n\right)_{n \in \mathbb{N}} \subset \mathscr{C}$ is uniformly Lipschitzian on $\Lambda$. According to the proof of Proposition $\ref{P34}$, for every $n \in \mathbb{N}$, and $0 \leqslant t_1<t_2 \leqslant T$, we have,
$$
\begin{aligned}
\left|\varphi_n\left(t_2\right)-\varphi_n\left(t_1\right)\right| & =\left|\mathbf{I}_{0^+}^{\alpha }\left(^{C}{\mathfrak{D}_{0^+}^{\alpha}} \varphi_n\left(t_2\right)\right)-\mathbf{I}_{0^+}^{\alpha }\left(^{C}{\mathfrak{D}_{0^+}^\alpha} \varphi_n\left(t_1\right)\right)\right| \\
& \leqslant \frac{3C\left(\overline{G}(1)\right)^\frac{1}{g^+}}{\Gamma(\alpha+1)}\left\|^{C}{\mathfrak{D}_{0^+}^{\alpha}} \varphi_n \right\|_{G}(t_2-t_1)^{\frac{\alpha}{g^+}}\\
&\leqslant \frac{3MC\left(\overline{G}(1)\right)^\frac{1}{g^+}}{\Gamma(\alpha+1)}(t_2-t_1)^{\frac{\alpha}{g^+}}.
\end{aligned}
$$
Therefore, by the Arzela-Ascoli theorem, $\left(\varphi_n\right)_{n \in \mathbb{N}}$ is relatively compact in $\mathscr{C}$. Consequently, there exists a subsequence of $\left(\varphi_n\right)_{n \in \mathbb{N}}$ converging strongly in $\mathscr{C}$ to $\varphi$, by the uniqueness of the weak limit. The proof is now complete.
\hfill$\Box$
\section{An application to  fractional boundary value problem }\label{sec4}
In this section, we study a fractional boundary problem involving the fractional Caputo-derivative operator in the new fractional Orlicz space $\mathcal{O}$. Specifically, we consider the following problem

\begin{equation*}\label{P}
	(\mathcal{P})~~\left\{\begin{array}{l}
^{C}{\mathfrak{D}_{T^-}^{\alpha}}\left(g\left(^{C}{\mathfrak{D}_{0^+}^{\alpha}} v(t) \right) \right)=a(t, v(t)),\quad \text{in }\Lambda, \\
v(0)=v(T)=0,
\end{array}\right.
\end{equation*}
 where $0<\alpha<1$, $g$ is as given in Section $\ref{sec1}$ and the right hand-side $a:\Lambda \times \mathbb{R} \rightarrow \mathbb{R}$ satisfies the following conditions :
 \begin{enumerate}
     \item[$(a_1$)]:  $a \in \mathcal{C}(\Lambda \times \mathbb{R})$.
     \item[$(a_2$)]: \label{a2}  There is a constant $\mu>k$, where $k$ is given by the $\Delta_2-$condition such that $0<\mu A(t, v) \leq v a(t, v)$ for every $t \in\Lambda$ and $v \in \mathbb{R} \backslash\{0\}$.
 \end{enumerate}
 The assumption \hyperref[a2]{$(a_2)$} is known by the Ambrossetti-Rabinowitz condition which is often a key assumption in the application of the mountain pass theorem, which is a fundamental tool in critical point theory for finding saddle points of functionals. 
Net, a function $v \in \mathcal{O}_{0}$ is a weak solution for \hyperref[P]{($\mathcal{P}$)} if
\begin{equation}
     \int_{\Lambda}g\left({}^{C}{\mathfrak{D}_{0^+}^{\alpha}} v(t) \right) {}^{C}{\mathfrak{D}_{0^+}^{\alpha}}\phi(t) d t=\int_{\Lambda} a(t, v(t)) \phi(t) d t
\end{equation}
for any $\phi\in \mathcal{O}_{0}$.\\
In order to investigate the existence of nontrivial week solution of  problem \hyperref[P]{($\mathcal{P}$)}, we consider the following energy functional 
\begin{equation}\label{42}
    \mathcal{J}(v)=\int_{\Lambda}G\left({}^{C}{\mathfrak{D}_{0^+}^{\alpha}} v(t) \right) d t - \int_{\Lambda} A(t, v(t))  d t,
    \end{equation}
where $A(t, x)=\int_0^x a(t, \xi) d \xi$.
\begin{lemma}
    The function $\mathcal{J} \in \mathcal{C}^1\left(\mathcal{O}_{0}, \mathbb{R}\right)$ and
for all $v, \phi \in \mathcal{O}_{0}$ $$
\left\langle \mathcal{J}^{\prime}(v), \phi\right\rangle=\int_{\Lambda}g\left({}^{C}{\mathfrak{D}_{0^+}^{\alpha}} v(t) \right) {}^{C}{\mathfrak{D}_{0^+}^{\alpha}}\phi(t) d t - \int_{\Lambda} a(t, v(t)) \phi(t) d t,
$$
 where $\left\langle \cdot , \cdot \right\rangle $ denotes the usual duality between $ \mathcal{O}_{0}$ and its dual space $ \mathcal{O}_{0}^*$.
\end{lemma}
\noindent \textbf{\textit{Proof.}} Let us define
$$\mathcal{I}(v)=\int_{\Lambda}G\left({}^{C}{\mathfrak{D}_{0^+}^{\alpha}} v(t) \right)   d t
\quad \text{ and } \quad
 \mathcal{E}(v)=\int_{\Lambda} A(t, v(t)).$$\\
First, it is easy to see that
$$
\left\langle \mathcal{J}^{\prime}(v), \phi\right\rangle=\int_{\Lambda}g\left({}^{C}{\mathfrak{D}_{0^+}^{\alpha}} v(t) \right) {}^{C}{\mathfrak{D}_{0   ^+}^{\alpha}}\phi(t) d t - \int_{\Lambda} a(t, v(t)) \phi(t) d t.
$$
Next, we prove that $\mathcal{I} \in \mathcal{C}^1\left(\mathcal{O}_{0}, \mathbb{R}\right)$.\\
Let $\left\{v_k\right\} \subset \mathcal{O}_0$ with $v_k \longrightarrow v$ strongly in $\mathcal{O}_0$, show that $\mathcal{I}^{\prime}\left(v_k\right) \longrightarrow \mathcal{I}^{\prime}(v)$ in $\mathcal{O}_0^*$.
Indeed,
$$
\left\langle \mathcal{I}^{\prime}(v_k)-\mathcal{I}^{\prime}(v), \phi\right\rangle=\int_{\Lambda}\left[g\left({}^{C}{\mathfrak{D}_{0^+}^{\alpha}} v_k(t) \right)-g\left({}^{C}{\mathfrak{D}_{0^+}^{\alpha}} v(t) \right)\right] {}^{C}{\mathfrak{D}_{0   ^+}^{\alpha}}\phi(t) d t .
$$
Let us set
$$
P_k=g\left({}^{C}{\mathfrak{D}_{0^+}^{\alpha}} v_k\right) \in L^{\overline{G}}(\Lambda),
$$
$$
P=g\left({}^{C}{\mathfrak{D}_{0^+}^{\alpha}} v\right) \in L^{\overline{G}}(\Lambda),
$$
$$
\overline{P}={}^{C}{\mathfrak{D}_{0   ^+}^{\alpha}}\phi(t) \in L^{G}(\Lambda).
$$
Hence, by the Hölder inequality and $(\ref{23:})$, we obtain
$$
\left\langle \mathcal{I}^{\prime}\left(v_k\right)-\mathcal{I}^{\prime}(v), \varphi\right\rangle \leqslant 2\left\|P_k-P\right\|_{L^{\overline{G}}(\Lambda)}\|{\overline{P}}\|_{L^{G}(\Lambda)}
$$
Thus
$$
\left\|\mathcal{I}^{\prime}\left(v_k\right)-\mathcal{I}^{\prime}(v)\right\|_{\mathcal{O}_0^*} \leqslant 2\left\|P_k-P\right\|_{L^{\overline{G}}(\Lambda)}.
$$
Now since $P_k \longrightarrow P$ in $\mathcal{O}_0$. Then $P_k \longrightarrow P$ in $L^{\overline{G}}(\Lambda)$.\\
Hence, for a subsequence of $\left(P_k\right)_{k \geqslant 0}$, we get $P_k(x) \longrightarrow P(x)$ a.e. in $\Lambda$ and $\exists h \in L^{G}(\Lambda)$ such that $\left|P_k(x)\right| \leqslant h(x)$.\\
So we have
$P_k(x) \longrightarrow P(x)$ a.e. in $\Lambda$ and $\left|P_k(x)\right|\leq \left|g(|h(x)|)\right|$
. Consequently
$$
\mathcal{I}^{\prime}\left(v_k\right) \longrightarrow \mathcal{I}^{\prime}(v) \text { in } \mathcal{O}_0^* \text {. }
$$
In the same way, we show that $\mathcal{E}^{\prime}\left(v_k\right) \longrightarrow \mathcal{E}^{\prime}(v).$\\
Then, ih the light of the  Dominated Convergence Theorem, we deduce the continuity of $\mathcal{J}^{\prime}$.\hfill$\Box$\\

\noindent  Therefore, the week solution to \hyperref[P]{($\mathcal{P}$)} is determined by the critical points of $\mathcal{J}$.\\
Now we are ready to state our main existence result  as follows : 
\begin{theorem}\label{T41}
    Let $\alpha\in (0,1)$ and suppose that $a$ satisfy $(a_1)$ and $(a_2)$. The problem \hyperref[P]{($\mathcal{P}$)} has a nontrivial weak solution $v \in  \mathcal{O}_{0}.$
\end{theorem}
To establish the proof of Theorem \ref{T41}, we apply the mountain-pass theorem.\\
Now, we  recall and prove some auxiliaries results which are crucial in the proof of the main existence result of this section.
\begin{lemma}[\cite{Torres}]\label{L42}
    If $a$ satisfies $\left(a_2\right)$, then for every $t \in\Lambda$, the following inequalities hold
\begin{equation}\label{43}
A(t, v) \leq A\left(t, \frac{v}{|v|}\right)|v|^\mu, \text { if } 0<|v| \leq 1 ;
\end{equation}
and
\begin{equation}\label{44}
A(t, v) \geq A\left(t, \frac{v}{|v|}\right)|v|^\mu, \text { if }|v| \geq 1 .
\end{equation}

%%We note, according to Lemma $\ref{L42}$ $ f$ is superquadratic at infinity and subquadratic at the origin.
\end{lemma}
\begin{lemma}[\cite{Torres}]\label{L43}
     Let $\ell=\inf \{A(t, v) ~~|~~ t \in\Lambda,|v|=1\}$. Then for any $\xi \in \mathbb{R} \backslash\{0\}$ and $v \in \mathcal{O}_0$, we have
$$
\int_{\Lambda} A(t, \xi v(t)) d t \geq \ell|\xi|^\mu \int_{\Lambda}|v(t)|^\mu-T \ell.
$$
\end{lemma}
The following result known by the Palais-Smale compactness condition (PS) which  is a key assumption in the mountain-pass theorem.
\begin{lemma}
Let $G$ be an N-function which satisfies $(\ref{delta})$, we suppose that the function $a$ satisfies the conditions $(a_1)$ and $(a_2)$. Then the functional $\mathcal{J}$ given by $(\ref{42})$ satisfies the Palais-Smale condition.
\end{lemma}
\noindent \textbf{\textit{Proof.}} Let $\left\{v_k\right\}$ be a (PS)-sequence of $\mathcal{J}$ on $\mathcal{O}_{0}$, which could be expressed mathematically as follows
\begin{equation}\label{45}
    \left|\mathcal{J}\left(v_k\right)\right| \leq M\quad \text{and}\quad \lim _{k \rightarrow \infty} \mathcal{J}^{\prime}\left(v_k\right)=0 .
\end{equation}
We shall first demonstrate the boundedness of $\left\{v_k\right\}$. Note that
$$ \mathcal{J}(v_k)=\int_{\Lambda}G\left({}^{C}{\mathfrak{D}_{0^+}^{\alpha}} v_k(t) \right) d t - \int_{\Lambda} A(t, v_k(t))  d t,$$
and 
$$
 \left\langle \mathcal{J}^{\prime}(v_k), v_k\right\rangle =\int_{\Lambda}g\left({}^{C}{\mathfrak{D}_{0^+}^{\alpha}} v_k(t) \right) {}^{C}{\mathfrak{D}_{0   ^+}^{\alpha}}v_k(t) d t - \int_{\Lambda} a(t, v_k(t)) v_k(t) d t.
$$
Then by $(\ref{45})$, we get
\begin{equation}
    \begin{aligned}\left|\mathcal{J}\left(v_k\right)-\tfrac{1}{\mu} \left\langle \mathcal{J}^{\prime}\left(v_k\right), v_k\right\rangle \right| & \leq\left|\mathcal{J}\left(v_k\right)\right|+\left|\tfrac{1}{\mu} \mathcal{J}^{\prime}\left(v_k\right)\right|\left|v_k\right| \\ & \leq C\left(1+\left\|v_k\right\|_{\mathcal{O}_0}\right) .\end{aligned}
\end{equation}
However, with Proposition $\ref{P32}$ and $(\ref{Gg})$, we get

\begin{equation}\label{47}
    \begin{aligned} 
\mathcal{J}( v_k) & = \int_{\Lambda}G\left({}^{C}{\mathfrak{D}_{0^+}^{\alpha}} v_k(t) \right) d t - \int_{\Lambda} A(t, v_k(t))  d t \\ 
& \geqslant \left\|v_k\right\|_{\mathcal{O}_0}^{g^\pm}- \int_{\Lambda} A(t, v_k(t))  d t.
\end{aligned}
\end{equation}
In addition

\begin{equation}\label{48}
    \begin{aligned} 
   \left\langle  \mathcal{J}^{\prime}(v_k), v_k\right\rangle & = \int_{\Lambda}g\left({}^{C}{\mathfrak{D}_{0^+}^{\alpha}} v_k(t) \right) {}^{C}{\mathfrak{D}_{0   ^+}^{\alpha}}v_k(t) d t - \int_{\Lambda} a(t, v_k(t)) v_k(t) d t \\ & \leqslant k \int_{\Lambda} {}^{C}{\mathfrak{D}_{0   ^+}^{\alpha}}v_k(t) d t - \int_{\Lambda} a(t, v_k(t)) v_k(t) d t\\
    &\leqslant k \left\|v_k\right\|_{\mathcal{O}_0}^{g^\pm} - \int_{\Lambda} a(t, v_k(t)) v_k(t) d t.
    \end{aligned}
    \end{equation}
By $(\ref{47})$, $(\ref{48})$ and $(a_2)$, we have
\begin{equation}
    \begin{aligned}  \mathcal{J}\left(v_k\right)-\tfrac{1}{\mu} \left\langle \mathcal{J}^{\prime}\left(v_k\right), v_k\right\rangle   \geqslant & \left(1-\tfrac{k}{\mu}\right)\left\|v_k\right\|_{\mathcal{O}_0}^{g^\pm}-\int_{\Lambda} A\left(t, v_k(t)\right) d t+\tfrac{1}{\mu} \int_{\Lambda} a\left(t, v_k(t)\right) v_k(t) d t \\
    \geq & \left(1-\tfrac{k}{\mu}\right)\left\|v_k\right\|_{\mathcal{O}_0}^{g^\pm} .
    \end{aligned}
    \end{equation}
Since $\mu >k$ it follows that ${v_k}$ is bounded in $\mathcal{O}_0$. Since $\mathcal{O}_0$ is reflexive space, there exists some $v$ in $\mathcal{O}_0$ and some sub-sequence which for simplicity we still denote by $v_k$, with  $v_k \rightharpoonup v$ in $\mathcal{O}_0$, as $k \rightarrow \infty$.

\begin{equation}\label{410}
    \begin{aligned} 
\left\langle\mathcal{J}^{\prime}\left(v_k\right)-\mathcal{J}^{\prime}(v), v_k-v\right\rangle  = & \left\langle\mathcal{J}^{\prime}\left(v_k\right), v_k-v\right\rangle-\left\langle \mathcal{J}^{\prime}(v), v_k-v\right\rangle \\ \leq & \left\|\mathcal{J}^{\prime}\left(v_k\right)\right\|\left\|v_k-v\right\|_{\mathcal{O}_0}-\left\langle\mathcal{J}^{\prime}(v), v_k-v\right\rangle .
\end{aligned}
\end{equation}
By $(\ref{410})$, we find 
$$\left\langle\mathcal{J}^{\prime}\left(v_k\right)-\mathcal{J}^{\prime}(v), v_k-v\right\rangle \rightarrow 0, \text{ as } k \rightarrow \infty.
$$
From Propositions $\ref{P34}$ and $\ref{P35}$, we get that $v_k$ is bounded in $\mathscr{C}(\Lambda)$, additionally, we might presume that
$$\lim _{k \rightarrow \infty}\left\|v_k-v\right\|_{\infty}=0.$$
Hence, we have
$$
\int_{\Lambda}\left[a\left(t, v_k(t)\right)-a(t, v(t))\right]\left(v_k(t)-v(t)\right) d t \underset{k \rightarrow \infty}{\longrightarrow}0.
$$
Furthermore, a simple calculation reveals that 
$$
\left\langle \mathcal{J}^{\prime}\left(v_k\right)-\mathcal{J}^{\prime}(v), v_k-v\right\rangle\geqslant\left\|v_k-v\right\|_{\mathcal{O}_0}-\int_{\Lambda}\left(a\left(t, v_k(t)\right)-a(t, v(t))\right)\left(v_k(t)-v(t)\right) d t .
$$
So $\left\|v_k-v\right\|_{\mathcal{O}_0} \underset{k \rightarrow \infty}{\longrightarrow}0$. That is $\left\{v_k\right\}$ converges strongly to $v$ in $\mathcal{O}_0$.\hfill$\Box$\\

\noindent Now, we are in position to establish the proof of the main result of this section.\\

 \noindent \textbf{\textit{Proof of Theorem}} $\ref{T41}$. It is simple to observe that $\mathcal{J}(0)=0.$ It remains to show that $\mathcal{J}$ satisfies the geometric conditions of the mountain pass theorem. \\
 From $(\ref{IN})$, we infer $$\max _{t \in\Lambda}|v(t)| \leq M\|v\|_{\mathcal{O}_0}, \forall v \in \mathcal{O}_0,$$
with $M=\tfrac{C\left(\overline{G}(1)\right)^\frac{1}{g+}}{\Gamma(\alpha+1)}x^{\frac{\alpha}{g}}$. Next let $C_1=\tfrac{1}{C}$, it follows by the inequality from above
and $(\ref{43})$, if $\|v\|_{\mathcal{O}_0} \leq C_1$
$$
\begin{aligned}
\int_{\Lambda} A(t, v(t)) d t & \leq \int_{\Lambda} A\left(t, \frac{v(t)}{|v(t)|}\right)|v(t)|^\mu d t \\
& \leq MTC^{\mu}\|v\|_{\mathcal{O}_0}^{\mu}.
\end{aligned}
$$
Then
$$\begin{aligned}  \mathcal{J}(u)&=\int_{\Lambda}G\left({}^{C}{\mathfrak{D}_{0^+}^{\alpha}} v(t) \right) d t - \int_{\Lambda} A(t, v(t))  d t
\\ & \geq \|v\|_{\mathcal{O}_0}^{g^\pm}-MTC^{\mu}\|v\|_{\mathcal{O}_0}^{\mu}, \quad \text { if }\|v\|_{\mathcal{O}_0} \leq C_1,\end{aligned}$$
Therefore
$$\mathcal{J}(v)\geq C_1^{g^{\pm}}-MTC_1^{\mu}C^{\mu}.$$
Let us consider $R<\min \left\{C_1,\left(\tfrac{1}{ M T C^\mu}\right)^{\frac{1}{\mu-g^{\pm}}}\right\}$ and $\beta=R^{g^\pm}-M T C^\mu R^\mu$, then
$$
\mathcal{J}(v) \geq \beta ~~\text { with }\|v\|_{\mathcal{O}_0}=R.
$$
Thus $\mathcal{J}$ fulfils the first geometric condition of the mountain pass theorem.\\

\noindent Now, for each $\xi \in \mathbb{R} \backslash\{0\} $ and $v \in \mathcal{O}_0 $, from Lemma $\ref{L43}$, we infer that  
$$
\begin{aligned}
\mathcal{J}(\xi v) & \leqslant \xi \|v\|_{\mathcal{O}_0}^{g^\pm} - \int_{\Lambda} A(t, \xi v(t))  d t \\
& \leqslant \xi \|v\|_{\mathcal{O}_0}^{g^\pm}-\ell|\xi|^\mu \int_{\Lambda}|v(t)|^\mu d t+T \ell .
\end{aligned}
$$
 Since $\mu>k$, by passing to the limit when $\xi \rightarrow+\infty$, we obtain $\mathcal{J}(\xi v) \rightarrow-\infty$, so the second geometric condition of the mountain pass theorem is affirmed by taking $e=\xi v$, with $\xi$ sufficiently large, then $\mathcal{J}(e) \leq 0.$\\
 
 \noindent Therefore $\mathcal{J}$ satisfies the
mountain pass condition, and thus  $\mathcal{J}$ possess a nontrivial critical point which is a nontrivial weak solution of \hyperref[P]{($\mathcal{P}$)}.\\

%--------------------------------
\noindent\textbf{{\large Declarations}} :\\

\noindent\textbf{Ethical Approval :} Not applicable.\\

\noindent\textbf{Competing interests :} The authors declare that there is no conflict of interest.\\

\noindent\textbf{Authors' contributions :} The authors contributed equally to this work.\\

\noindent\textbf{Funding :} Not applicable.\\

\noindent\textbf{Availability of data and materials :} Not applicable.

    % % % % % % % % % % % % % % % % % % % % % % % % % % % % % % % % % % % % % % % % % % % % % % % % % % % % % % % % % % % % % % % % % % % % % % % % % % % % % % % % % % % % % % % % % % % % % % % %

\end{document}